\documentclass[11pt,graybox,envcountchap,oribibl]{svmono}

\usepackage{emptypage}
\usepackage[a5,pdflatex,center]{crop}
\usepackage{mathptmx}
\usepackage{helvet}
\usepackage{courier}
\usepackage{type1cm}
\DeclareMathAlphabet{\mathcal}{OMS}{cmsy}{m}{n}

\usepackage{makeidx}
\usepackage{multicol}
\usepackage{graphicx}
\usepackage[bottom]{footmisc} 

\usepackage[english]{babel}
\usepackage{amssymb,amsmath,enumitem,interval,graphbox,url,cancel,float,pdfpages}
\usepackage{tikz}
\usetikzlibrary{decorations.fractals,decorations.markings,arrows.meta,calc,intersections,shapes.geometric,decorations.pathreplacing,angles,quotes}

\newlength{\raggio} 
\newlength{\mra}
\newlength{\bmra}
\newlength{\golden}
\colorlet{TA}{cyan!15}
\colorlet{TB}{green!15}
\colorlet{SA}{magenta!10}
\colorlet{SB}{yellow!15}
\colorlet{GTA}{gray!10}
\colorlet{GTB}{gray!60}
\colorlet{GSA}{gray!10}
\colorlet{GSB}{gray!60}

\newcommand{\N}{\mathbb{N}}
\newcommand{\Z}{\mathbb{Z}}

\newcommand{\C}{\mathbb{C}}

\usepackage{etoolbox}
\AtEndEnvironment{proof}{\hfill\ensuremath{\blacksquare}}

\usepackage{mathtools}
\mathtoolsset{showonlyrefs} 

\makeindex

\makeatletter
\def\@maketitle{\newpage
 \null
 \vskip 2em                 
\begingroup
  \def\and{\unskip, }
  \parindent=\z@
  \pretolerance=10000
  \rightskip=\z@ \@plus 3cm
  {\LARGE                   
   \lineskip .5em
   \@author
   \par}%
  \vskip 1cm                
  {\Huge \@title \par}
  \vskip 3mm                
  \if!\@subtitle!\else
   {\LARGE\ignorespaces\@subtitle \par}
   \vskip 1cm                
  \fi
	\vfill
	\begin{center}
	\includegraphics[scale=0.8]{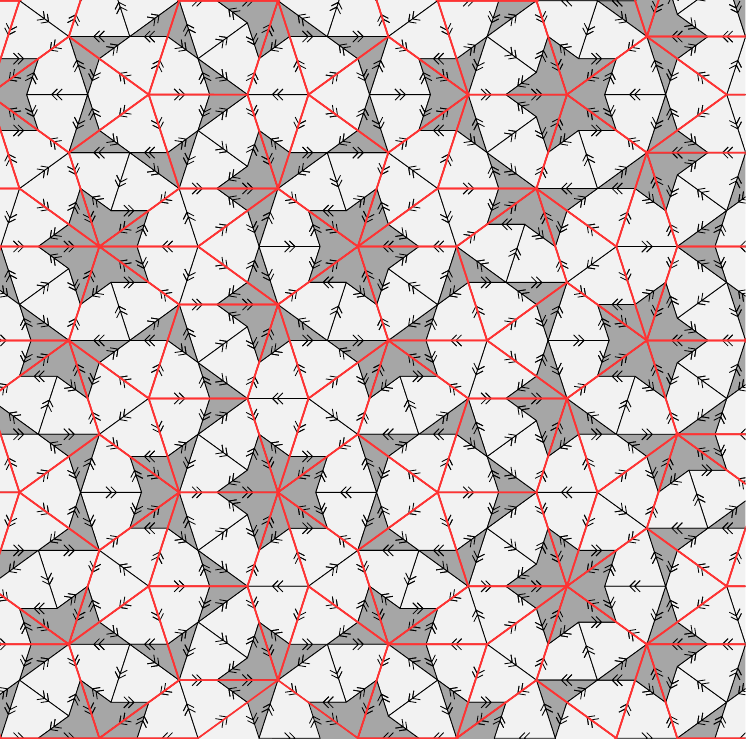}
	\end{center}
	\vfill
  \if!\@date!\else
    {\large \@date}
    \par
    \vskip 1.5em               
  \fi
\endgroup}
\makeatother

\begin{document}

\author{Francesco D'Andrea}
\title{A guide to Penrose tilings} 
\date{April 24, 2023}

\begin{titlepage}
\maketitle
\end{titlepage}

\frontmatter

\preface

The aim of this book is to provide an elementary introduction, complete with detailed proofs, to the celebrated tilings of the plane discovered by Sir Roger Penrose in the `70s. The book covers many aspects of Penrose tilings, including the study of the space parameterizing Penrose tilings from the point of view of Connes' Noncommutative Geometry.

I am in debt with my students and colleagues for their comments on a preliminary version of this text. A special thanks goes to Prof.\ Giovanni Landi for his suggestions about the last chapter of the book.

All the images in this book were created in \LaTeX\ with the TikZ package.

\vspace{2cm}

\begin{itshape}
This is a preview of the following work: Francesco D'Andrea, ``A Guide to Penrose Tilings'', 2023, Springer, ISBN: 978-3-031-28427-4.

It is reproduced with permission of Springer.
\end{itshape}

\cleardoublepage

\includepdf[width=148mm,pages=-]{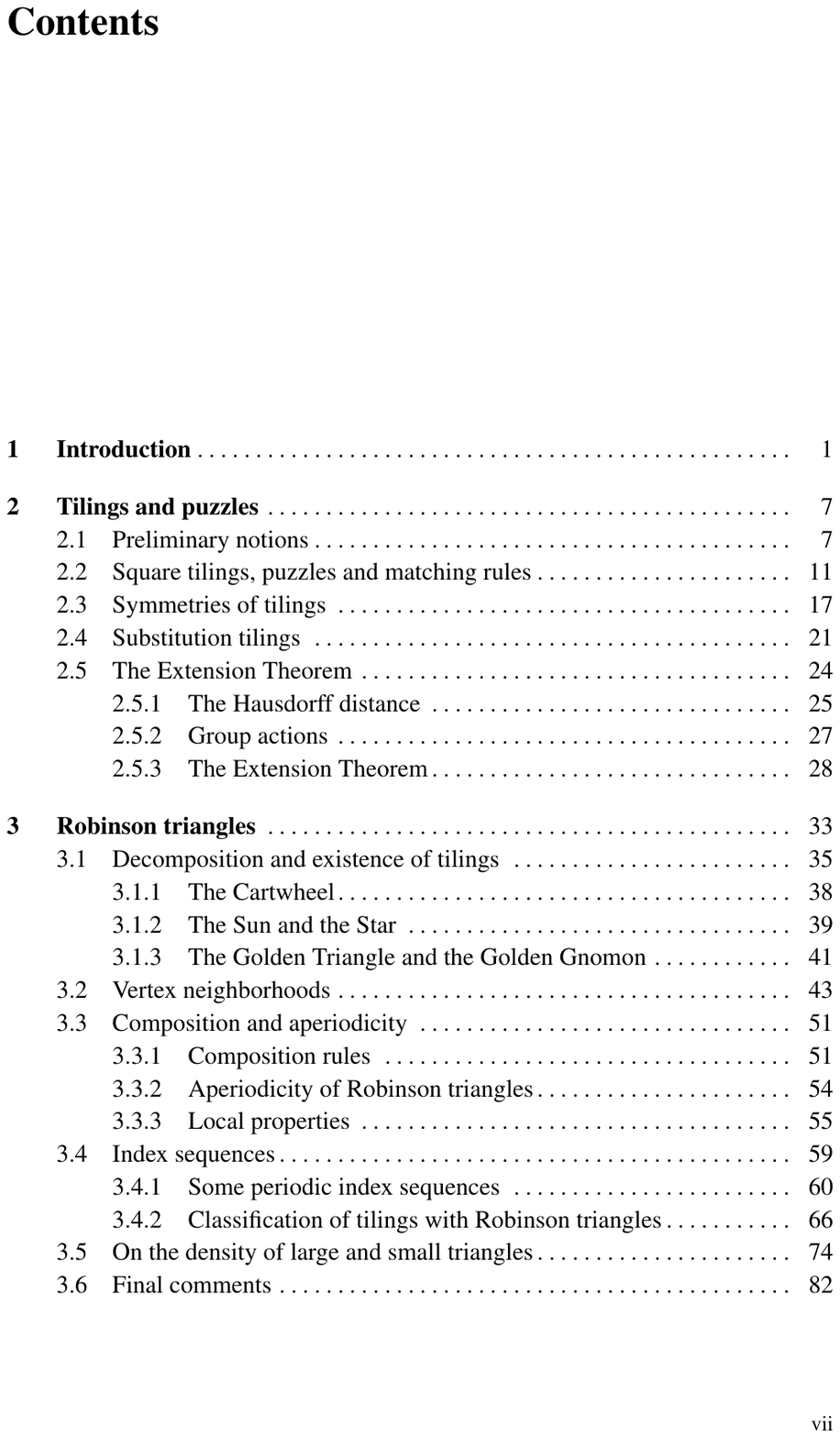}

\mainmatter

\chapter{Introduction}

The problem of covering a flat surface --- a subset of the Euclidean plane or the whole plane itself ---
using some fixed geometric shapes and with no overlaps is probably one of the oldest in mathematics.
Such a covering is called a \emph{tiling}, or also a \emph{tessellation} or a \emph{mosaic} (see Def.~2.1).
The picture below shows a tiling with equilateral triangles and (non-regular) pentagons.

\begin{figure}[h!]
\begin{center}
\includegraphics{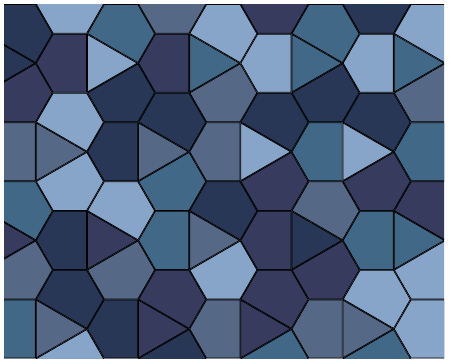}
\end{center}
\caption{A portion of a tiling of the plane.}\label{fig:mosaic}
\end{figure}

Mosaics with colored geometrical shapes can be seen in stained glass windows in Christian churches. Tilings using shapes that are invariant under a rotation of $72^\circ$, such as regular pentagons or pentagrams, occur frequently in Islamic art.
These shapes are said to possess a 5-fold symmetry (since $72$ is one-fifth of $360$).

\smallskip

Probably, one of the reasons for the popularity of this topic is that many of its fundamental problems are simple to formulate, even if they may require advanced tools to be solved. Not always, though:
the story of the mathematical amateur Marjorie Rice, illustrated e.g.~in \cite{Sc78}, is a famous example of
a non-professional mathematician giving an important original contribution to the topic.

\smallskip

Despite the apparently random placement of tiles in Figure \ref{fig:mosaic}, an optical effect due to the fancy choice of colors, it is not difficult to realize that one can extend the tiling to the whole plane simply by translating copies of the same basic configuration
\begin{center}
\includegraphics[page=1]{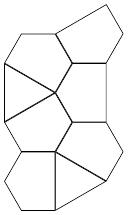}
\end{center}
In this way, one gets a first (not the easiest) example of \emph{periodic} tiling, where the same pattern is repeated over and over again.

\smallskip

In 1961, Wang proposed the following problem \cite{Wan61}. Consider a finite number of square tiles, with edges marked by symbols (or colors, letters, numbers, \ldots). We want to tile the whole plane by placing copies of these tiles, which we call \emph{Wang tiles}, in the plane. The tiles can be translated but not rotated and two adjacent tiles must share a full edge: this is what we call an \emph{edge-to-edge} tiling. Moreover, there is a \emph{matching rule}: two edges can touch only if they are decorated with the same symbol. The Domino Problem asks if there exists an algorithm that can decide whether a set of Wang tiles can tile the whole plane. Wang proved that the Domino Problem is decidable if and only if there does \underline{not} exist an \emph{aperiodic} set of Wang tiles, i.e.~a set that can tile the plane only non-periodically, and for this reason he conjectured that such a set didn't exist. The first aperiodic set of Wang tiles was found by Berger in `66 \cite{Ber66} and consisted of more than 20 thousand different tiles (but Berger's original Ph.D.\ Thesis contained a simplified version with 104 tiles). A notable contribution to this topic is the one by Knuth, the creator of \TeX, who found an aperiodic set of 92 Wang tiles \cite[Sect.~2.3.4.3]{Knu68}. Many others contributed to the subject (we will not discuss the Domino Problem in this book, but the interested reader can find more about it in \cite{GS87}).
The number of Wang tiles in an aperiodic set was recently reduced to 11, and it was proved that this is the minimal number \cite{JR15}.

\smallskip

The story of the Domino Problem intertwines with another one: the search for tilings of the plane that use only shapes with 5-fold symmetry. It is well known that this cannot be achieved with a single tile \cite{DGS82}, for example a regular pentagon, or by repeating the same pattern over and over again (see Sect.~2.3). This problem can be found already in Kepler's
treatise \textit{Harmonices mundi} (1619), where one can find several mosaics with regular pentagons, decagons, and pentagrams.
Of course, one can tile the whole plane with a single pentagonal tile if the angles are chosen correctly. For example, this can be achieved with the pentagons in Figure \ref{fig:mosaic} by repeating the pattern
\begin{center}
\includegraphics[page=2]{mosaic3.pdf}
\end{center}
over and over again. Kepler proved that, up to a similarity,
there exist exactly eleven different edge-to-edge tilings of the plane with convex regular polygons that are semiregular, i.e.~that have only one type of vertex neighborhood (see Sect.~2.1 for the definition). Three of them are made with a single polygon (i.e.~they are \emph{monohedral}): either an equilateral triangle, a square, or a regular hexagon.

\smallskip

In 1973, looking for tilings of the plane with shapes with 5-fold symmetry, Penrose found a set of 6 tiles that can only tile the plane non-periodically \cite{Pen74}, thus reducing the cardinality of known aperiodic sets to 6. One year later, he reduced this number to 2, by discovering his famous tiles shaped like a \emph{kite} and a \emph{dart} (see Chap.~4 for the details):
\begin{center}

\begin{tikzpicture}[semithick]

\coordinate (v1) at (0,0);
\coordinate (v2) at (54:1.618);
\coordinate (v3) at (0,1.618);
\coordinate (v4) at (126:1.618);
\coordinate (v5) at (4,1.4);
\coordinate (v6) at ($(v5)+(-54:1.618)$);
\coordinate (v7) at ($(v5)+(0,-1)$);
\coordinate (v8) at ($(v5)+(234:1.618)$);
\coordinate (v9) at ($(v5)+(0,-2.5)$);

\draw[gray!50,thick,smooth,tension=0.6] plot coordinates{(v1) (-0.2,-0.4) (0.2,-0.8) (0,-1.2)};

\draw[gray!50,thick] (v7) -- (v9);

\foreach \k in {0,...,3} \draw[thick,gray!50] ($(v9)+(0,0.15*\k)$) -- ++(-30:0.2) ($(v9)+(0,0.15*\k)$) -- ++(210:0.2);

\draw (v1) -- (v2) -- (v3) -- (v4) -- cycle;
\draw (v5) -- (v6) -- (v7) -- (v8) -- cycle;
\draw[dashed,gray] (v1) -- (v3) (v5) -- (v7);


\end{tikzpicture}
\end{center}
Even if the discovery was in 1973-74, the first official mention was only in 1977 in \cite{Gar77} (as explained therein, Penrose applied for patents in the UK, US and Japan, and was reluctant to advertise the tilings before the patents were granted). The paper \cite{Pen78} contains an account of how Penrose discovered his sets with 6 and 2 tiles,
and a nice Escheresque non-periodic tiling with birds (now called ``Penrose non-periodic chickens''): it was produced using the same technique employed by Escher to transform polygonal tiles into animal shapes.

Let us mention that Socolar and Taylor recently found an aperiodic disconnected monotile \cite{ST11,ST12}.
Even more recently, in the preprint \cite{SMKGS23}, the Authors exhibit a tile homeomorphic to a disk that, when used together with its reflection, can tile the plane only non-periodically (the interesting story of this discovery can be found, e.g., in \cite{Qua23}).
An overview of previously known aperiodic sets of tiles can be found in Chap.~7 of \cite{Sen95}.

There exists many books on Penrose tilings intended for a general audience (e.g.~\cite{Gar97}), survey papers (e.g.~the original papers by Penrose \cite{Pen74,Pen78}) and even some excellent 
YouTube videos (e.g.~\cite{Veritasium}) carefully explaining all the properties of these tilings.

An efficient way to study Penrose tilings is by dissecting them into triangles:
\begin{equation}\label{eq:dissection}
\begin{tikzpicture}[semithick,baseline=(current bounding box.center)]

\coordinate (v1) at (0,0);
\coordinate (v2) at (54:1.618);
\coordinate (v3) at (0,1.618);
\coordinate (v4) at (126:1.618);
\coordinate (v5) at (4,1.4);
\coordinate (v6) at ($(v5)+(-54:1.618)$);
\coordinate (v7) at ($(v5)+(0,-1)$);
\coordinate (v8) at ($(v5)+(234:1.618)$);

\draw (v1) -- (v2) -- (v3);
\draw[transform canvas={xshift=-6pt}] (v3) -- (v4) -- (v1);
\draw (v5) -- (v6) -- (v7);
\draw[transform canvas={xshift=-6pt}] (v7) -- (v8) -- (v5);
\draw[dashed,gray] (v1) -- (v3) (v5) -- (v7);
\draw[dashed,gray,transform canvas={xshift=-6pt}] (v1) -- (v3) (v5) -- (v7);

\end{tikzpicture}
\end{equation}
This was the method used by Robinson in \cite{Rob75}, that we illustrate in Chap.~3. Robinson triangles are extremely useful to study local properties of Penrose tilings. Moreover, they allow to introduce a fundamental tool in the classification of Penrose tilings up to isometries: the \emph{index sequence} of a tiling (Sect.~3.4).

A few years later, another fundamental tool for the construction of Penrose tilings of the plane was discovered by de Bruijn: the so-called \emph{pentagrid} method \cite{dBr81}.
A very nice historical account of the impact of de Bruijn on this subject is in \cite{Yang13}.
Rather than tiling the plane by constructing bigger and bigger patches, de Bruijn discovered a global construction using five bundles of parallel lines. A byproduct of this method is the interpretation of Penrose non-periodic tilings as two-dimensional slices of periodic tilings of a five-dimensional space. At the core of de Bruijn's idea is the observation that the projection
$$
\Z^n\to\C,\qquad
(k_0,\ldots,k_{n-1})\mapsto \sum_{j=0}^{n-1}k_j(e^{2\pi\mathrm{i}/n})^j ,
\qquad(n\geq 1)
$$
of an $n$-dimensional lattice on a plane is a lattice if $n\in\{1,2,3,4,6\}$, and has dense image in $\C$ in all other cases. For $n=5$, if we project only points that are ``near''
the plane, we get a rhombus tiling that is equivalent (in a suitable sense) to a Penrose tiling by kites and darts. The image of the above projection is the ring of cyclotomic integers $\Z[e^{2\pi\mathrm{i}/n}]$, and the pentagrid method is a bridge between Penrose tilings and algebraic number theory.

An important stimulus for the study of non-periodic tilings came from a discovery in chemistry.
In 1982, Dan Shechtman observed a diffraction pattern with ten-fold symmetry 
from a metal alloy sample. Such a symmetry is impossible in ordinary crystals, that are modeled on lattices (cf.~Sect.~2.3).
For his work, Shechtman was awarded the Nobel Prize in Chemistry in 2011 for the discovery of quasi-crystals \cite{Fer18}.

A (doubly) periodic tiling is formed by repeating a bounded region infinitely many times. It follows that every finite subset of tiles (every ``patch'' in the terminology that we will introduce in the next chapter) is repeated infinitely many times in the tiling. A tiling with the latter property is called \emph{repetitive}, and surprisingly a tiling does not need to be periodic to be repetitive. Penrose tilings, for example, are non-periodic and repetitive. A tiling that is non-periodic and repetitive is called \emph{quasi-periodic}.
Quasi-periodic tilings have attracted the interest of many scientists because they provide a mathematical model for quasi-crystals.

\smallskip

Using index sequences one shows that the space parametrizing inequivalent Penrose tilings is the quotient of a Cantor space by an equivalence relation with dense equivalence classes. The quotient topology is then the trivial (indiscrete) one, and gives no insight on the structure of the space of Penrose tilings.
In the spirit of Noncommutative Geometry \cite{Con94}, quotient spaces that are not Hausdorff can be studied using groupoid C*-algebras, and the space of equivalence classes of Penrose tilings is an example of such a ``bad'' quotient space.
It is therefore natural to use C*-algebra techniques to learn more about these tilings.
The use of C*-algebra techniques in solid state physics was popularized by Bellissard in the `80s \cite{Bel86}. A good starting point for their use in the study of tilings is the beautiful survey by Kellendonk and Putnam \cite{KP00} (see also \cite{Sad08}).

\bigskip

This book is structured as follows. The second chapter is a general introduction to tilings and illustrates the main notions through some simple examples.
The third chapter is about Robinson's triangles. The fourth chapter is about Penrose tilings.
The fifth chapter is about de Bruijn's pentagrid method, and the realization of a Penrose tiling as a projection of a portion of a five-dimensional lattice.
The sixth and last chapter is about the noncommutative geometry of Penrose tilings.
The geometric properties Robinson triangles and Penrose tiles are intimately related to the algebraic properties of the \emph{golden ratio}. Some of these algebraic properties are collected in Appendices A.1 and A.2.

\bigskip

We shall adopt the following notations. We denote by $\N$ the set of natural numbers including $0$.
If $S$ is a set, $|S|$ denotes its cardinality. If $S$ is a subset of a topological space, $\overline{S}$ is its closure and $\mathring{S}$ its interior. $A\subseteq B$ or $B\supseteq A$ means that $A$ is a subset of $B$, while $A\subset B$ or $B\supset A$ means that $A$ is a proper subset of $B$.
An open interval with endpoints $a$ and $b$, is denoted by $\ointerval{a}{b}$, a closed one by $\interval{a}{b}$, etc.

\cleardoublepage

\thispagestyle{empty}

\phantom{.}\vfill

\begin{center}
PAGES FROM 7 TO 198 ARE NOT INCLUDED IN THIS PREVIEW
\end{center}

\vfill\phantom{.}

\backmatter

\setcounter{page}{199}

\end{document}